\documentclass[10pt,journal,onecolumn]{IEEEtran}

\usepackage{graphicx, bm, amsmath, siunitx, longtable,tabularx, float, mathtools, booktabs, eqnarray, color, amsfonts, amssymb}
\usepackage[version=4]{mhchem}
\newcommand{\B}[1]{{\bm #1}}
\newcommand{\T}{^{\mbox{\tiny T}}}
\newcommand{\ds}{\displaystyle}
\newcommand{\dd}{\; \text{d}}
\newcommand{\bdot}[1]{\dot{\B{#1}}}

\newcommand{\tfc}{Theory of Functional Connections}
\newcommand{\ces}{constrained expressions}
\renewcommand{\ce}{constrained expression}
\newcommand{\p}[2]{\prescript{(#1)}{}{#2}}

\usepackage{etoolbox}

\begin{document}

\title{Fuel-Efficient Powered Descent Guidance on Large Planetary Bodies via Theory of Functional Connections$^*$ \thanks{$^*$This work was partially supported by a NASA Space Technology Research Fellowship, Johnston [NSTRF 2019] Grant \#: 80NSSC19K1149.}
}

\author{Hunter Johnston$^1$\thanks{$^1$Aerospace Engineering, Texas A\&M University, College Station, TX, 77843 USA (hunterjohnston@tamu.edu, mortari@tamu.edu)} \and Enrico Schiassi$^2$\thanks{$^2$Systems and Industrial Engineering, The University of Arizona (eschiassi@email.arizona.edu)} \and Roberto Furfaro$^{2,3}$\thanks{$^3$Aerospace and Mechanical Engineering, The University of Arizona (robertof@email.arizona.edu)} \and Daniele Mortari$^1$}

\maketitle

\begin{abstract}
    In this paper we present a new approach to solve the fuel-efficient powered descent guidance problem on large planetary bodies with no atmosphere (e.g. the Moon or Mars) using the recently developed Theory of Functional Connections. The problem is formulated using the indirect method which casts the optimal guidance problem as a system of nonlinear two-point boundary value problems. Using the Theory of Functional Connections, the problem constraints are analytically embedded into a ``constrained expression,'' which maintains a free-function that is expanded using orthogonal polynomials with unknown coefficients. The constraints are satisfied regardless of the values of the unknown coefficients which convert the two-point boundary value problem into an unconstrained optimization problem. This process casts the solution into the admissible subspace of the problem and therefore simple numerical techniques can be used (i.e. in this paper a nonlinear least-squares method is used). In addition to the derivation of this technique, the method is validated in two scenarios and the results are compared to those obtained by the general purpose optimal control software, GPOPS-II. In general, the proposed technique produces solutions of $\mathcal{O}(10^{-10})$. Additionally, for the proposed test cases, it is reported that each individual TFC-based inner-loop iteration converges within 6 iterations, each iteration exhibiting a computational time between 72 and 81 milliseconds within the MATLAB legacy implementation. Consequently, the proposed methodology is potentially suitable for on-board generation of optimal trajectories in real-time.
\end{abstract}

\section{Introduction}
Precision landing on large planetary bodies is a technology of utmost importance for future human and robotic exploration of the solar system. Over the past two decades, landing systems for robotic Mars missions have been developed and successfully deployed robotic assets on the Martian surface (e.g. rovers \cite{grotzinger2012mars}, landers \cite{shotwell2005phoenix}). Considering the strong interest in sending humans to Mars within the next few decades, as well as the renewed interest in building infrastructure in the Earth-Moon system for easy access to the Lunar surface, the landing system technology will need to progress to satisfy the  demand for more stringent requirements. One of the enabling technology for precision landing is the ability to generate on-board and track in real-time fuel optimal trajectories. Generally, finding optimal guidance solutions requires setting up a cost function (i.e. the objective of the optimization) and the desired control and state constraints. For landing on large planetary bodies, especially for systems that deliver large systems to their surface, it is extremely important to minimize the fuel. State and control constraints usually include the system dynamics (i.e. equations of motion), thrust value limitation as well as state space constraints that need to be enforced to ensure safety (e.g., glide slope constraints to ensure  safety of the spacecraft). 

Usually, two methods are available to solve optimal control problems, direct and indirect methods. Direct methods are based on discretizing the continuous states and controls with the goal of transforming the continuous problem into a Non Linear Programming (NLP) problem \cite{darby2011hp,fahroo2002direct,ross2004pseudospectral}. The latter can be cast as a finite constrained optimization problem that can be solved via any of the available numerical algorithms that have the potential to find a local minimum (e.g., thrust region method \cite{byrd2000trust}). Whereas direct methods have been applied to solve a large variety of optimal control problems (e.g. \cite{josselyn2003rapid,graham2015minimum,miller2017rapid,jiang2019integrated}), the general NLP problem is considered to be NP-hard, i.e. non deterministic polynomial time hard. NP-hard problems imply that the required amount of computational time needed to find the optimal solution does not have a predetermined bound (i.e. a bound cannot be apriori determined). NP-hard problems are such that the computational time necessary to converge to the solution is not known. As a consequence, the lack of assured convergence may result in questioning the reliability of the proposed approach. Since for optimal closed loop space guidance, most of the problems require computing numerical solutions on board and in real time, general algorithms that solve NLP problems cannot be reliably implemented. More recently, researchers have been experimenting with transforming optimal control problem from a general non convex formulation into a convex optimization problem \cite{acikmese2007convex,blackmore2010minimum}. Here, the goal is to take advantage of the assured convex convergence properties. Indeed, convex optimization problems are shown to be computationally tractable as their related numerical algorithms guarantee convergence to a global optimal solution in a polynomial time. The general convex methodology requires that the optimal guidance problem is formulated as convex optimization whenever appropriate or convexification techniques are applied to transform the problem from a non-convex into a convex one. Such methodologies have been proposed and applied to solve optimal guidance and control via direct method in a large variety of problems including, planetary landing \cite{acikmese2007convex,blackmore2010minimum}, entry atmospheric guidance \cite{wang2016constrained,wang2018autonomous}, rocket ascent guidance \cite{zhang2019rapid}, and low thrust \cite{wang2018minimum}.

Alternatively, a second approach to solve optimal control and guidance problems has been generally applied to a variety of optimal control problems. Named indirect method, the approach applies optimal control theory (i.e. Pontryagin Minimum Principle) to formally derive the first-order necessary conditions that must be satisfied by the optimal solution (state and control). The problem is cast as a Two Point Boundary Value Problem (TPBVP) that must be solved to determine the time evolution of state and costate from which the control generally depends. For general nonlinear problems, the necessary conditions result in a complicate set of equations and conditions. The resulting TPBVP tends to be highly sensitive to the initial guess on the costates making the problem very hard to solve. Although indirect methods are known to yield more accurate optimal solutions, they are very hard to implement and tend to be less used in practice (with respect to direct methods).

Recently, a new approach \cite{U-TFC,M-TFC} called \emph{Theory of Functional Connections} (TFC)\footnote{This theory, initially called ``Theory of Connections'', has been renamed for two reasons. First, the ``Theory of Connections'' already identifies a specific theory in differential geometry, and second, what this theory is actually doing is ``Functional Interpolation'' as it provides \emph{all} functions satisfying a set of constraints in term of function and any derivative in rectangular domains of $n$-dimensional space.} was developed to derive expressions, called \emph{constrained expression}, with embedded constraints. This approach has been successfully applied to solve both linear \cite{LDE} and nonlinear \cite{NDE} differential equations, and IVP and BVP, at machine error accuracy and in milliseconds. TFC is a general methodology providing \emph{functional interpolation} expressions with an embedded set of $n$ linear constraints.  Such expressions can be expressed in the following general form:
\begin{equation}\label{eq01}
    y (t) = g (t) + \ds\sum_{k=1}^n \eta_k \, p_k (t) 
\end{equation}
Here, the $p_k (t)$ are $n$ assigned linearly independent functions and the $g (t)$ is a free function that must be linearly independent from the selected $p_k (t)$ terms. The $\eta_k$ are coefficient functions that are derived by imposing the set of $n$ constraints. The constraints considered in the TFC are any linear combination of the functions and/or derivatives evaluated at specified values of the variable $t$. The constraints of the differential equation to be solved enable to directly derive the unknown coefficient functions, $\eta_k$. Once the $\eta_k$ terms are determined, the constraints for the differential equations are satisfied for any possible $g (t)$. 

The constrained expressions thus obtained can be used to transform the constrained optimization problem in an unconstrained one. This implies reducing the whole solution search space to just the space of admissible solutions  (i.e. those fully complying with all constraints). In the field of optimization, this has been done by expanding the free function by a set of basis functions (e.g., Fourier series or orthogonal polynomials, such as Legendre or Chebyshev polynomials) whose coefficients are found by direct application of a least-squares algorithms. Nonlinear initial or boundary value problems require the implementation of an iterative least-squares approach to converge to the desired solution \cite{NDE}.

In this paper, we propose solving the optimal landing guidance problem via a new method based on TFC that can compute fuel-efficient trajectories both fastly and accurately. After deriving the TPBVP arising from the fuel-efficient powered descent guidance necessary conditions, the TFC is employed to generate the whole set of boundary condition-free equations that can be solved by expanding the solution in Chebyshev polynomials and computing the expansion coefficients using Iterative Least-Squares (ILS) method. The proposed methodology is shown to be fast, accurate and thus potentially suitable for on-board generation of optimal landing trajectories on large planetary bodies.

This article is organized in the following manner. First, a summary of the TFC approach to solve a general TPBVP is presented. Next, the equations for the optimal powered descent landing problem are presented and the necessary conditions are derived. Afterward, it is shown how the TFC approach is used to transform the TPBVP and eventually solved via and iterative least-squares approach. Lastly, results for two specific optimal thrust profiles are presented: test \#1) min-max profile and test \#2) max-min-max profile. These results are compared against those obtained from the general purpose optimal control software GPOPS-II \cite{GPOPS}.

\section{The \tfc\ approach to solve a general TPBVP}\label{sec:tfc_example}
The uni-variate TFC \cite{U-TFC} embeds $n$ linear constraints into the expression,
\begin{equation*}
    y(t) = g(t) + \sum_{k=1}^n \eta_k \, p_k(t) 
\end{equation*}
where $g(t)$ is a free function, $\eta_k$ are unknown coefficients to be solved, and $p_k(t)$ are user defined linearly independent functions; where $g(t)$ and the set of $p_k(t)$ must be linearly independent. This formulation is easily extended for the use in solving vector equations, by utilizing index notation where 
\begin{equation}\label{eq:CE_vector}
    y_i(t) = g_i(t) + \sum_{k=1}^n \eta_{k_i} \, p_k(t) \quad \text{for} \quad i = 1,2,3,
\end{equation}
such that $i$ denotes each component of the vector. Since landing problems are second-order BVPs the following example derivation will follow that form. For this class of problems, we have the conditions $y_i(t_0) = y_{0_i}$, $\dot{y}_i(t_0) = \dot{y}_{0_i}$, $y_i(t_f) = y_{f_i}$, and $\dot{y}_i(t_f) = \dot{y}_{f_i}$. For these constraint conditions we can select $p_k (t) = t^{(k-1)}$ (see Ref. \cite{LDE,U-TFC,NDE} for detailed explanation on this choice) in Eq. (\ref{eq:CE_vector}), we obtain,
\begin{equation*}
    \begin{aligned}
        \begin{Bmatrix} y_{0_i}-g_{0_i} \\ y_{f_i}-g_{f_i} \\ \dot{y}_{0_i}-\dot{g}_{0_i} \\ \dot{y}_{f_i}-\dot{g}_{f_i}\end{Bmatrix} = \begin{bmatrix} 1 & t_0 & t_0^2 & t_0^3 \\ 1 & t_f & t_f^2 & t_f^3 \\ 0 & 1 & 2t_0 & 3t_0^2 \\ 0 & 1 & 2 t_f & 3t_f^2\end{bmatrix} \begin{Bmatrix} \eta_{1_i} \\ \eta_{2_i} \\ \eta_{3_i} \\ \eta_{4_i}\end{Bmatrix}
    \end{aligned}
\end{equation*}
which can be solved for the $\eta_{k_i}$ terms through matrix inversion which leads to the \ce\ of the form,
\begin{equation}\label{eq:vect_tpbvp}
\begin{aligned}
    y_i(t) = g_i(t) &+ \Omega_1(t) \Big(y_{0_i} - g_{0_i}\Big) + \Omega_2(t) \Big(y_{f_i} - g_{f_i}\Big) + \\ &+ \Omega_3(t) \Big(\dot{y}_{0_i} - \dot{g}_{0_i}\Big) + \Omega_4(t) \Big(\dot{y}_{f_i} - \dot{g}_{f_i}\Big)
\end{aligned}
\end{equation}
where the $\Omega_k(t)$ terms are defined as in Tables \ref{tab:om_general} where we define $t_* = t - t_0$ and $\Delta t = t_f - t_0$
\begin{table}[H]
\begin{center}
\begin{tabular}{SSSSS} \toprule
{} & {Initial Value} & {Final Value} & {Initial Derivative} & {Final Derivative}\\
{} & {$\Omega_1(t)$} & {$\Omega_2(t)$} & {$\Omega_3(t)$} & {$\Omega_4(t)$} \\ \bottomrule \midrule
{$(\cdot)$} & {$1 + \dfrac{2t_*^3}{\Delta t^3} - \dfrac{3 t_*^2}{\Delta t^2}$} & {$- \dfrac{2 t_*^3}{\Delta t^3} + \dfrac{3 t_*^2}{\Delta t^2}$} & {$t_* + \dfrac{t_*^3}{\Delta t^2} - \dfrac{2 t_*^2}{\Delta t}$} & {$\dfrac{t_*^3}{\Delta t^2} - \dfrac{t_*^2}{\Delta t}$} \\ \midrule
{$\dfrac{\dd}{\dd t}(\cdot)$} & {$\dfrac{6t_*^2}{\Delta t^3} - \dfrac{6t_*}{\Delta t^2}$} & {$-\dfrac{6t_*^2}{\Delta t^3} + \dfrac{6t_*}{\Delta t^2}$} & {$1 + \dfrac{3t_*^2}{\Delta t^2} - \dfrac{4t_*}{\Delta t}$} & {$\dfrac{3t_*^2}{\Delta t^2} - \dfrac{2t_*}{\Delta t}$} \\ \midrule
{$\dfrac{\dd^2}{\dd t^2}(\cdot)$} & {$\dfrac{12t_*}{\Delta t^3} - \dfrac{6}{\Delta t^2}$} & {$-\dfrac{12t_*}{\Delta t^3} + \dfrac{6}{\Delta t^2}$} & {$\dfrac{6t_*}{\Delta t^2} - \dfrac{4}{\Delta t}$} & {$\dfrac{6t_*}{\Delta t^2} - \dfrac{2}{\Delta t}$}\\ \midrule
\bottomrule
\end{tabular}
\end{center}
\caption{Switching functions for TPBVP constraints for general domain of $t \in [t_0, t_f]$.}
\label{tab:om_general}
\end{table}
Eq. (\ref{eq:vect_tpbvp}) represents all possible functions satisfying the boundary value constraints. Furthermore, the derivatives follow,
\begin{equation*}
    \begin{cases}
    \dot{y}_i(t) = \dot{g}_i(t) &+ \dot{\Omega}_1(t) \Big(y_{0_i} - g_{0_i}\Big) + \dot{\Omega}_2(t) \Big(y_{f_i} - g_{f_i}\Big) + \\ &+ \dot{\Omega}_3(t) \Big(\dot{y}_{0_i} - \dot{g}_{0_i}\Big) + \dot{\Omega}_4(t) \Big(\dot{y}_{f_i} - \dot{g}_{f_i}\Big) \\
    \ddot{y}_i(t) = \ddot{g}_i(t) &+ \ddot{\Omega}_1(t) \Big(y_{0_i} - g_{0_i}\Big) + \ddot{\Omega}_2(t) \Big(y_{f_i} - g_{f_i}\Big) + \\  &+ \ddot{\Omega}_3(t) \Big(\dot{y}_{0_i} - \dot{g}_{0_i}\Big) + \ddot{\Omega}_4(t) \Big(\dot{y}_{f_i} - \dot{g}_{f_i}\Big)
    \end{cases}
\end{equation*}
The \ce\ defined by Eq. (\ref{eq:vect_tpbvp}) and its derivatives can then be applied to a differential equation which, in general, can be expressed as a loss function (i.e. in its implicit form) such that,
\begin{equation}\label{eq:gen:DE}
    F_i\left(t,y_j,\dot{y}_j,\ddot{y}_j\right) = 0 \quad \text{for} \quad i,j = 1,2,3
\end{equation}
where each component of the differential equation can be a function of all component $j$. By substituting Eq. (\ref{eq:vect_tpbvp}) into Eq. (\ref{eq:gen:DE}), the differential equation is transformed to a new differential equation that we define as $\tilde{F}$, which is only a function of the independent variable $t$, the free-function, and its derivatives where,
\begin{equation*}
    \tilde{F}_i\left(t,g_j,\dot{g}_j,\ddot{g}_j\right) = 0 \quad \text{for} \quad i,j = 1,2,3.
\end{equation*}
This differential equation is unique\footnote{This differential equation, which has no constants of integration, has the boundary conditions embedded through the free function evaluated at the constraints.} because it is subject to \emph{no constraints} and  will always satisfy the boundary-values. In order to solve this problem numerically, we expand the free function $g_i(t)$ in terms of some known basis ($\B{h} (z)$) with unknown coefficients ($\B{\xi}$) such that,
\begin{equation}\label{eq:basis}
    g_i (t) = \B{\xi}_i \T\B{h} (z) \quad \text{where} \quad z = z(t)
\end{equation}
where $\B{\xi}_i$ are $m \times 1$ vectors of unknown coefficients where $i$ is the dimension and $m$ is the number of basis functions. In general the basis functions are defined on a different domains (Chebyshev and Legendre polynomials are defined on $z\in[-1,+1]$, Fourier series is defined on $z\in[-\pi,+\pi]$, etc.) so these functions must be linearly mapped to the independent variable $t$. This can be done using the equations, 
\begin{equation*}
z = z_0 + \frac{z_f-z_0}{t_f-t_0}(t - t_0) \quad \longleftrightarrow \quad t = t_0 + \frac{t_f-t_0}{z_f-z_0}(z - z_0),
\end{equation*}
where $t_f$ represents the upper integration limit. The subsequent derivatives the the free-function defined in Eq. (\ref{eq:basis}) follow,
\begin{equation*}
    \frac{\dd^{n} g_i}{\dd t^{n}} = \B{\xi}_i \T  \frac{\dd^{n} \B{h}(z)}{\dd z^{n}} \left(\frac{\dd z}{\dd t}\right)^{n},
\end{equation*}
where by defining,
\begin{equation*}
c := \frac{\dd z}{\dd t} = \frac{z_f - z_0}{t_f - t_0}
\end{equation*}
the expression can be simplified to, 
\begin{equation}\label{eq:basis_derv}
    \frac{\dd^{n} g_i}{\dd t^{n}} = c^{n} \B{\xi}_i \T  \frac{\dd^{n} \B{h}(z)}{\dd z^{n}},
\end{equation}
which defines all mappings of the free-function. Lastly, the domain $t\in[t_0,t_f]$ must be discretized by $N$ points. In this paper (as in prior papers) we consider the linear basis $\B{h}(z)$ as Chebyshev orthogonal polynomials. The optimal distribution of points is provided by collocation points \cite{Colloc,ChebCol}, defined as,
\begin{equation}\label{eq:collo}
    z_k = -\cos\left(\frac{k \pi}{N}\right) \quad \text{for} \quad k = 1, 2, \cdots, N,
\end{equation}
As compared to the uniform distribution point, the collocation point distribution allows a much slower increase of the condition number as the number of basis functions, $m$, increases. By defining the free function in this way and then discretizing the domain of the differential equations, $\tilde{F}_i$ becomes
\begin{equation*}
    \tilde{F}_i(t, \B{\xi}_j) = 0 \qquad \text{for} \qquad i,j = 1,2,3,
\end{equation*}
which are some set of functions (linear or nonlinear) of the unknown parameters $\B{\xi}_j$ of which many unconstrained optimization schemes can be applied. For the convenience of the reader, the TFC method is summarized below in Fig. \ref{fig:tfc_summary}
\begin{figure}[H]
    \centering\includegraphics[width=.75\linewidth]{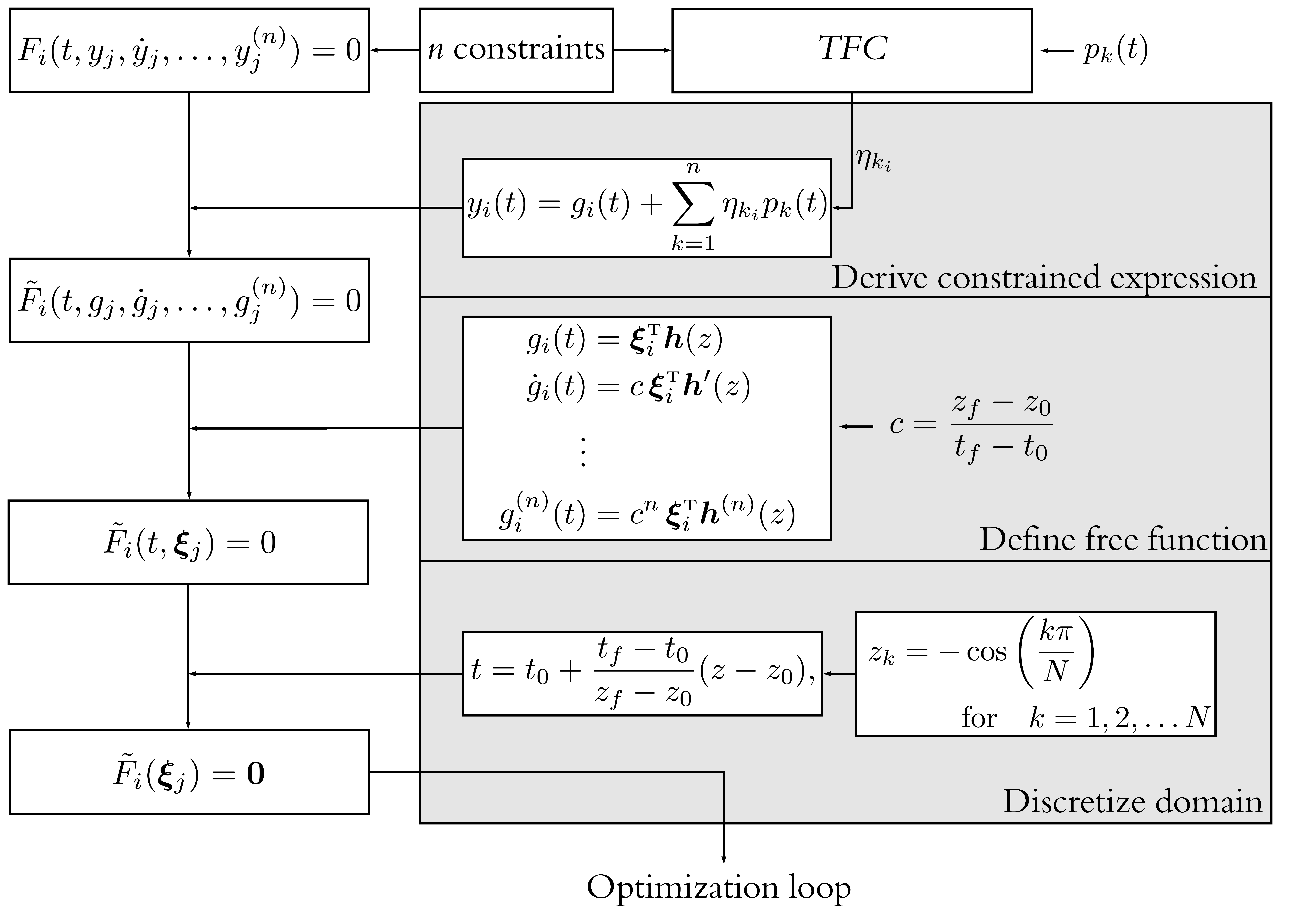}
    \caption{Summary of the relevant steps to use the TFC approach. The major steps of TFC include 1) deriving the constrained expression, 2) defining an appropriate free function, and 3) discretizing the domain.}
    \label{fig:tfc_summary}
\end{figure}

\section{Optimal Powered Descent Pinpoint Landing Problem}
For the problem of powered descent pinpoint landing guidance on large bodies (e.g. the Moon or Mars) the governing system dynamics during the powered descent phase can be modeled as:
\begin{eqnarray}
  \dot{\B{r}} &&= \,\B{v} \\ 
  \dot{\B{v}} &&= \,\B{a}_g + { \dfrac{\B{T}}{m} } \\ 
  \dot{m} &&= \,-\alpha \, T 
\end{eqnarray}
where $\alpha = 1 / v_{ex}$, $v_{ex}$ is the effective exhaust velocity of the rocket engine that is considered constant \cite{acikmese2007convex,upg}, $T = ||\B{T}||$, and $\B{T} = T \, \hat{\B{t}}$ is the thrust and it is constrained as follow:
\begin{equation*}
\begin{aligned}
    0 \leq T_{min} &\leq T \leq T_{max} \\
    || \hat{\B{t}} ||&=1
\end{aligned}
\end{equation*}
\begin{figure}[H]
    \centering\includegraphics[width=.75\linewidth]{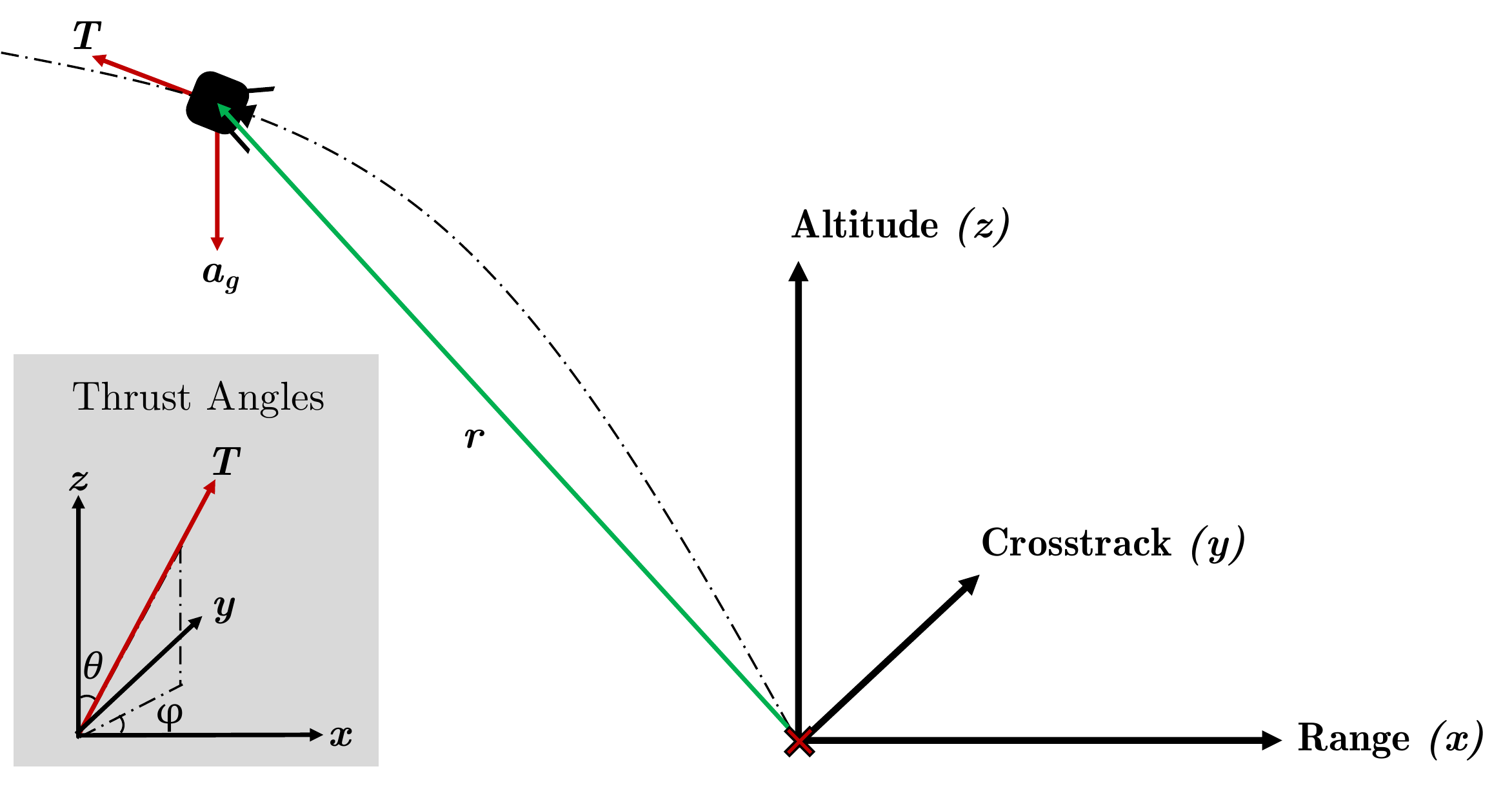}
    \caption{Coordinate frame definition for optimal powered descent pinpoint landing problem.}
    \label{fig:coord}
\end{figure}
Additionally, $\B{a}_g$ is the gravity acceleration which is also considered constant. As stated in \cite{upg}, this assumption is justified for short flights as it is the case for the powered descent phase of the landing process. A summary of the reference frame for the problem is given in Fig. \ref{fig:coord}. For this problem, the initial and final position and velocity, and initial mass are given:
\begin{equation*}
    \left\{\begin{split} \B{r}(0) = \B{r}_0 \\ \B{v}(0) = \B{v}_0\end{split}\right., \qquad \left\{\begin{split} \B{r}(t_f) = \B{r}_f \\ \B{v}(t_f) = \B{v}_f\end{split}\right., \qquad \text{and} \qquad m(0) = m_0.
\end{equation*}
The objective is to minimize the mass of propellant used while satisfying the dynamics constraints of the problem. Therefore the problem can be posed as,
\begin{equation*}
    \underset{t_f,T}{\text{minimize}} \quad \alpha\int_{0}^{t_f}T \dd \tau
\end{equation*}
\begin{equation*}
\begin{aligned}
\text{subject to} \quad \dot{\B{r}} = \,\B{v}, \quad &\dot{\B{v}} = \,\B{a}_g + { \dfrac{\B{T}}{m} }, \quad \dot{m} = \,-\alpha \, T \\
& 0 \leq T_{min} \leq T \leq T_{max}, \\
& \B{r}(0) = \B{r}_0, \quad \B{v}(0) = \B{v}_0, \quad m(0) = m_0 \\
& \B{r}(t_f) = \B{r}_f, \quad \B{v}(t_f) = \B{v}_f
\end{aligned}
\end{equation*}

\subsection{First-order necessary conditions (Pontryagin Minimum Principle)}
The application of the Pontryagin Minimum Principle (PMP) dictates that the Hamiltonian takes the following form \cite{Bryson_Ho},
\begin{equation*}
    H = L + \B{\lambda}\T \B{f} + \B{\mu}\T \B{C}
\end{equation*}
which can be expanded to,
\begin{equation}\label{eq:ham}
    H=\alpha T + \B{\lambda}_r\T \B{v}+ \B{\lambda}_v\T  \left( \B{a}_g + \frac{T}{m} \hat{\B{t}} \right) - \lambda_m \alpha T + \mu_1 (T-T_{\max}) + \mu_2 (T_{\min}-T)
\end{equation}
where $T-T_{\max}$ and $T_{\min}-T \leq 0$ and $\mu_1, \mu_2 > 0$. The control optimality is ensured as follows,
\begin{equation*}
    \frac{\partial H}{\partial \hat{\B{t}}} = 0 \quad \rightarrow \quad \text{minimum when} \; \B{\lambda}_v\T\hat{\B{t}} = -1 \quad \rightarrow \quad \hat{\B{t}} = -\frac{\B{\lambda}_v}{||\B{\lambda}_v||}
\end{equation*}
Thus Eq. (\ref{eq:ham}) can be rewritten as,
\begin{equation*}
H = \alpha T + \B{\lambda}_r\T \B{v}+ \B{\lambda}_v\T \B{a}_g - \frac{T}{m} ||\B{\lambda}_v|| - \lambda_m \alpha T + \mu_1 (T-T_{\max}) + \mu_2 (T_{\min}-T)
\end{equation*}
Now, in order to determine optimal thrust magnitude, we impose the following,
\begin{equation*}
    \frac{\partial H}{\partial T} = \alpha - \frac{1}{m}||\B{\lambda}_v|| - \alpha \lambda_m + \mu_1 - \mu_2 = 0 
\end{equation*}
There are three potential cases:
\begin{enumerate}
    \item if $\mu_1 = \mu_2 = 0 \quad (T_{\min} < T < T_{\max})$ \quad \text{then} \quad $\alpha - \frac{1}{m}||\B{\lambda}_v|| - \alpha \lambda_m = \sigma = 0$
    \vspace{1mm}
    \item if $\mu_1 = 0$, $\mu_2 > 0 \quad (T=T_{\min})$  \. \quad \qquad \text{then} \quad $\sigma - \mu_2 = 0 \rightarrow \sigma = \mu_2 > 0$
    \vspace{1mm}
    \item if $\mu_1 > 0$, $\mu_2 = 0 \quad (T = T_{\max})$ \: \: \quad \quad \text{then} \quad $\sigma + \mu_1 = 0 \rightarrow \sigma = -\mu_1 < 0$
\end{enumerate}
Finally, one can conclude that the thrust magnitude has the following program:
\begin{equation*}
    T = \begin{cases} = T_{\max} \qquad \text{if} \qquad \sigma < 0 \\ = T_{\min} \qquad \text{if} \qquad \sigma > 0 \end{cases}
\end{equation*}
It has been demonstrated \cite{upg} that the singular case $\sigma = 0$ corresponds to a constant thrust perpendicular to the gravity vector which is generally not possible for a powered descent problem. Therefore, a singular arc is not part of the sought optimal solution. Furthermore, it has been shown that for the guided powered descent on large planetary bodies, the switching function $\sigma$ changes signs at most twice \cite{upg}. Consequently, the thrust magnitude can switch between min-max twice at the most. That is, in the most general case, the thrust magnitude has a bang-bang profile max-min-max. Hence, we can write the thrust magnitude as function of time with $t_1$ and $t_2$ as parameters, where $t_1$ and $t_2$ are the times where the switches happen; i.e. $T=T(t;t_1,t_2)$.

The additional first-order necessary conditions for the costate equations are written as,
\begin{eqnarray*}
    \bdot{\lambda}_r &= -\dfrac{\partial H }{\partial \B{r}} &= \B{0} \\
    \bdot{\lambda}_v &= -\dfrac{\partial H }{\partial \B{v}} &= -\bdot{\lambda}_r \\
    \dot{\lambda}_m &= -\dfrac{\partial H }{\partial m} &= -\frac{T}{m^2}||\B{\lambda}_v||
\end{eqnarray*}
and the transversality condition implies that,
\begin{equation*}
    \lambda_m(t_f) = 0 \qquad \text{and} \qquad H(t_f) = 0.
\end{equation*}

\section{Solution via the Theory of Functional Connections}

With the simplifications introduced in the previous section, in order to find the optimal state and thrust program, the following non-linear TPBVP must be solved:
\begin{eqnarray}
  \bdot{r} &&= \B{v} \label{eq:r_dot}\\
  \bdot{v} &&= \B{a}_g - \beta(t)
  \dfrac{\B{\lambda_v}}{||\B{\lambda_v}||}\label{eq:v_dot}\\
  \bdot{\lambda}_r &&= \B{0} \label{eq:lam_r_dot}\\ 
  \bdot{\lambda}_v &&=  -\B{\lambda_r}\label{eq:lam_v_dot} \\
  \dot{\lambda}_m &&= \,- \frac{T (t;t_1,t_2)}{m^2}||\B{\lambda}_v|| \label{eq:lam_m_dot}\\
 H(t_f) =  0 &&= \alpha T(t_f;t_1,t_2) + \B{\lambda}_v\T \left(\B{a}_g - \beta(t_f)\frac{\B{\lambda}_v}{||\B{\lambda}_v||} \right) \label{eq:trans}
\end{eqnarray}
where $\beta(t) \triangleq \frac{T(t;t_1,t_2)}{m(t)}$ subject to $\B{r(0)} = \,\B{r_0} , \B{v(0)} = \,\B{v_0}$ and $\B{r(t_f)} = \,\B{r_f}, \B{v(t_f)} = \,\B{v_f}$. It must be noted that $\lambda_m$, does not show up in any other equation and therefore Eq. (\ref{eq:lam_m_dot}) can be solved independently. Since the the transversality condition gives $\lambda_m (t_f) = 0$, Eq. (\ref{eq:lam_m_dot}) can be solved by back propagation or by simply using the TFC method. The solution of this equation is provided in the Appendix.

Since the solution of this problem exhibits a bang-bang profile for thrust, the original formulation of the TFC method (i.e. that derived in Section \ref{sec:tfc_example}) must be adjusted to accommodate switching behavior in the control. The general theory for this extension to hybrid systems is studied in \cite{pieceTFC}, and it is applied in the following derivation. Additionally, there are a few equations that are redundant and can be removed completely via the TFC \ce\ to further simplify the solution of this nonlinear system of equations. First, TFC \ce\ are analytical expressions meaning that derivative of a \ce\ for $\B{r}(t)$ is exactly the function $\B{v}(t)$. Therefore, the differential equation expressed by Eq. (\ref{eq:r_dot}) is unnecessary and can be disregarded. Similarly, the equations for $\bdot{\lambda}_r$ and $\bdot{\lambda}_v$ can be simplified. First, let us express the vector equations as three scalar equations each where the index $i$ represents the individual components. Using this notation we can expand $\B{\lambda}_v$ such that,
\begin{equation}\label{eq:ce_costate}
    \lambda_{v_i}  = a_{0_i} + a_{1_i} z = \B{h}_{\lambda}\T \B{\xi}_{\lambda_i}, \quad \text{for} \quad i = 1,2,3
\end{equation}
which satisfies Eqs. (\ref{eq:lam_r_dot}-\ref{eq:lam_v_dot}) through
\begin{align*}
    \dot{\lambda}_{v_i} &= a_{1_i} \\
    \dot{\lambda}_{v_i} = -\lambda_{r_i} &= -a_{1_i} 
\end{align*}
This process reduces the problem to the solution of a single differential equation expressed by Eq. (\ref{eq:v_dot}) and an algebraic equation for the Hamiltonian at final time given by Eq. (\ref{eq:trans}). Rewriting the differential equation in index notation and collecting all terms on one side, a loss function based on the residuals of the differential equation can be defined,
\begin{equation}\label{eq:loss}
{\cal L}_i = a_i - a_{g_i} + \beta(t) \, \lambda_{v_i} \left(\displaystyle\sum_{j=1}^3 \lambda^2_{v_j}\right)^{-1/2} \qquad \text{for} \qquad i = 1,2,3
\end{equation}
Now, the only step left is to construct a \ce\ for the state variables. Importantly, care must be taken in the construction of such \ce. In the above derivation of the thrust structure, we have shown that the thrust switches at most twice leading to a max-min-max profile. Therefore, the function $\beta(t)$ in Eq. (\ref{eq:loss}) jumps twice along the solution trajectory. This switching causes three distinct differential equations that cannot be solved with a single polynomial expansion over the entire domain, as it is done for the energy-optimal guidance \cite{TFC-Control}. Therefore, a new formulation for the TFC approach has been developed to handle these hybrid systems \cite{pieceTFC}. This process allows for the continuity between each segment of the domain.
\begin{figure}[h]
    \centering\includegraphics[width=.65\linewidth]{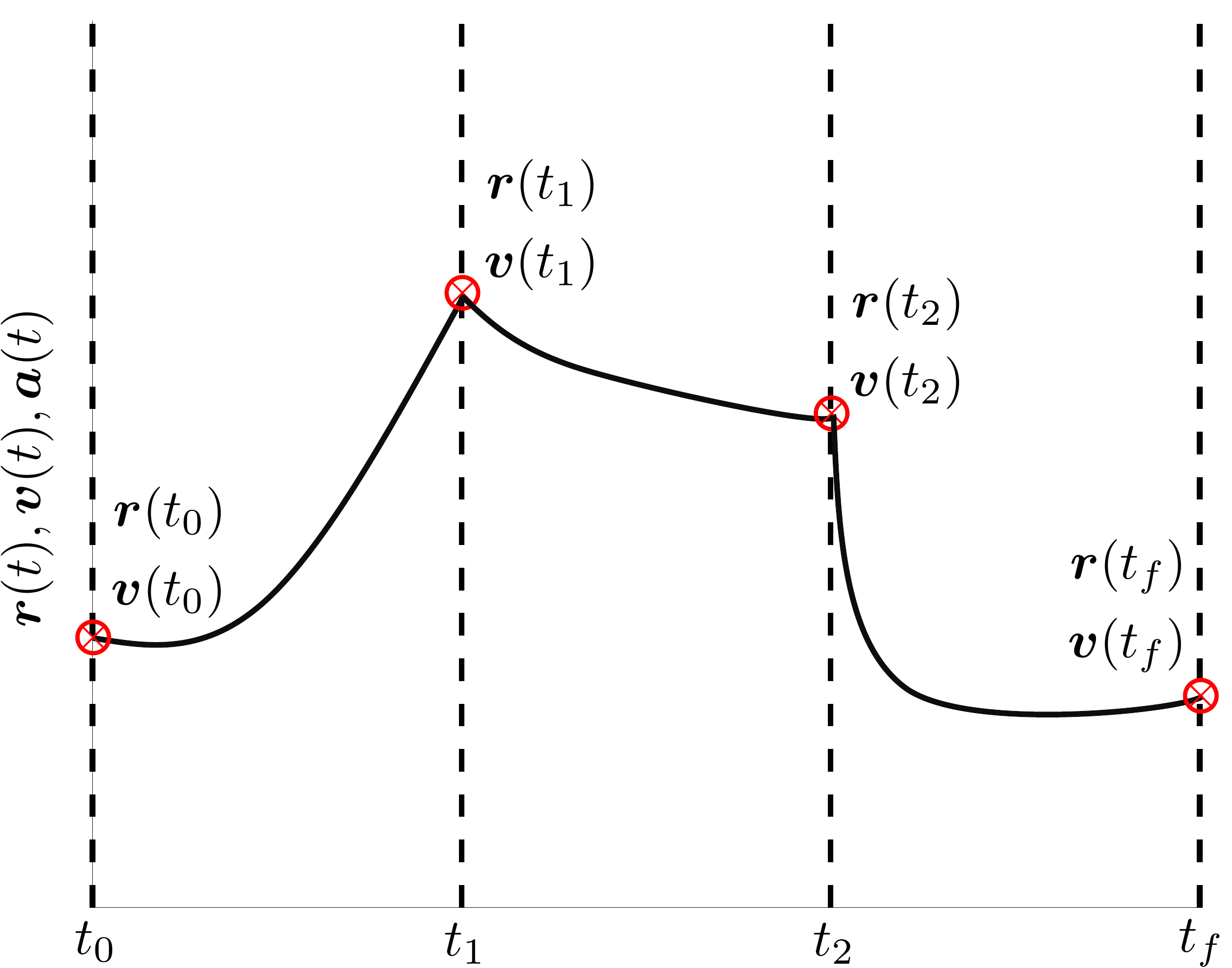}
    \caption{Visual representation of piece-wise approach using the TFC method \cite{pieceTFC}. In this derivation, the \ces\ maintain continuity of position and velocity through embedded relative constraints.}
    \label{fig:piece}
\end{figure}
As shown in Fig. \ref{fig:piece}, it is apparent that all sub-domains share the same constraint conditions (i.e, the initial and final position and velocity are constrained). Therefore, a single constraint expression can be derived for the case of arbitrary constraints and then incorporated into the sub-domains. The \ce\ for this specific case was derived in Section \ref{sec:tfc_example} and it is captured by Eq. (\ref{eq:vect_tpbvp}). Consequently, the position, velocity, and acceleration \ce\ can be expressed as,
\begin{eqnarray}
r_i = g_i \; + & \!\! \Omega_1(r_{0_i} - g_{0_i}) + \Omega_2 (r_{f_i} - g_{f_i}) + \nonumber \\ ~&+ \Omega_3 (v_{0_i} - \dot{g}_{0_i}) + \Omega_4(v_{f_i} - \dot{g}_{f_i})  \label{eq:ce_r} \\ v_i = \dot{g}_i\; + & \!\! \dot{\Omega}_1(r_{0_i} - g_{0_i}) + \dot{\Omega}_2 (r_{f_i} - g_{f_i}) + \nonumber \\ ~&+ \dot{\Omega}_3 (v_{0_i} - \dot{g}_{0_i}) + \dot{\Omega}_4(v_{f_i} -\dot{g}_{f_i}) \label{eq:ce_v} \\a_i = \ddot{g}_i \; + & \!\!  \ddot{\Omega}_1(r_{0_i} - g_{0_i}) + \ddot{\Omega}_2 (r_{f_i} - g_{f_i}) + \nonumber \\ ~&+ \ddot{\Omega}_3 (v_{0_i} - \dot{g}_{0_i}) + \ddot{\Omega}_4(v_{f_i} -\dot{g}_{f_i}) \label{eq:ce_a} 
\end{eqnarray}
where the $\Omega$ parameters are solely a function of the independent variable and act as switching functions to force the expression to always satisfy the specified constraints. These functions and their associated derivatives are summarized in Table \ref{tab:om_general} where $t_*$ and $\Delta t$ are defined by the specific segment. The \ce\ detailed by Eq. (\ref{eq:ce_r}-\ref{eq:ce_a}) can be used as a template to write the \ces\ for each segment of the solution trajectory. In order to explicitly identify the segment, the pre-superscript notation will be used. For example, $\p{1}{r}_i$ describes the position \ce\ for the first segment defined on $t\in[t_0,t_1]$. For this problem, $s=1$ is defined on $t\in[t_0,t_1]$, $s=2$ is defined on $t\in[t_1,t_2]$, and $s=3$ is defined on $t\in[t_2,t_f]$. Using this formulation, the \ces\ of position for each segment are,
\begin{equation*}
\begin{aligned}
\p{1}{r}_i = \p{1}{g}_i &+  \p{1}{\Omega}_1\left(r_{0_i} - \p{1}{g}_{0_i}\right) + \p{1}{\Omega}_2 \left(r_{1_i} - \p{1}{g}_{f_i}\right) + \\ &+ \p{1}{\Omega}_3 \left(v_{0_i} - \p{1}{\dot{g}}_{0_i}\right) + \p{1}{\Omega}_4\left(v_{1_i} -\p{1}{\dot{g}}_{f_i}\right) \\
\p{2}{r}_i = \p{2}{g}_i &+  \p{2}{\Omega}_1\left(r_{1_i} - \p{2}{g}_{0_i}\right) + \p{2}{\Omega}_2 \left(r_{2_i} - \p{2}{g}_{f_i}\right) + \\ &+ \p{2}{\Omega}_3 \left(v_{1_i} - \p{2}{\dot{g}}_{0_i}\right) + \p{2}{\Omega}_4\left(v_{2_i} -\p{2}{\dot{g}}_{f_i}\right) \\
\p{3}{r}_i = \p{3}{g}_i &+  \p{3}{\Omega}_1\left(r_{2_i} - \p{3}{g}_{0_i}\right) + \p{3}{\Omega}_2 \left(r_{f_i} - \p{3}{g}_{f_i}\right) + \\ &+ \p{3}{\Omega}_3 \left(v_{2_i} - \p{3}{\dot{g}}_{0_i}\right) + \p{3}{\Omega}_4\left(v_{f_i} -\p{3}{\dot{g}}_{f_i}\right)
\end{aligned}
\end{equation*}
where the derivative of these functions follow the form of Eq. (\ref{eq:ce_r}-\ref{eq:ce_a}). Now, the function of $g_i$ can also be expressed as a linear basis such that,
\begin{equation*}
    g_i(t) = \B{h}(z)\T \B{\xi}_i \quad \text{for} \quad i = 1,2,3.
\end{equation*}
with the subsequent derivatives following according to Eq. (\ref{eq:basis_derv}). This allows us to collect the unknown $\B{\xi}_i$ vectors and write the \ce\ in the form,
\begin{equation*}
\begin{aligned}
\p{1}{r}_i = \p{1}{} &\Big(\B{h} - \Omega_1\B{h}_0 - \Omega_2\B{h}_f - \Omega_3\bdot{h}_0 - \Omega_4\bdot{h}_f \Big)\T \p{1}{\B{\xi}}_i + \\ &+ \p{1}{\Omega}_1 r_{0_i} + \p{1}{\Omega}_2 r_{1_i} + \p{1}{\Omega}_3 v_{0_i} + \p{1}{\Omega}_4 v_{1_i} \\
\p{2}{r}_i = \p{2}{} &\Big(\B{h} - \Omega_1\B{h}_0 - \Omega_2\B{h}_f - \Omega_3\bdot{h}_0 - \Omega_4\bdot{h}_f \Big)\T \p{2}{\B{\xi}}_i + \\ &+ \p{2}{\Omega}_1 r_{1_i} + \p{2}{\Omega}_2 r_{2_i} + \p{2}{\Omega}_3 v_{1_i} + \p{2}{\Omega}_4 v_{2_i} \\
\p{3}{r}_i = \p{3}{} &\Big(\B{h} - \Omega_1\B{h}_0 - \Omega_2\B{h}_f - \Omega_3\bdot{h}_0 - \Omega_4\bdot{h}_f \Big)\T \p{3}{\B{\xi}}_i + \\ &+ \p{3}{\Omega}_1 r_{2_i} + \p{3}{\Omega}_2 r_{f_i} + \p{3}{\Omega}_3 v_{2_i} + \p{3}{\Omega}_4 v_{f_i} \\
\end{aligned}
\end{equation*}
Along with the linear unknowns in $\B{\xi}_i$, the equations share linear unknowns in $r_{1_i},v_{1_i},r_{2_i},v_{2_i}$ which serve as the embedded relative constraints between adjacent segments. With this new formulation, we now have three separate loss functions based on the residual of the differential equation over each segment ($s$) which are as follows,
\begin{equation*}\label{eq:loss_all_seg}
\p{s}{{\cal L}}_i = \p{s}{a}_i - a_{g_i} + \beta(t) \, \lambda_{v_i} \, \left(\displaystyle\sum_{j=1}^3 \lambda^2_{v_j}\right)^{-1/2}.
\end{equation*}
Note that although the costate \ces\ do not need to be split into separate domains,  special attention must be paid to discretizing the equations according to the segment time ranges. In order to solve for the unknown $\B{\xi}_i$ parameters, a nonlinear least-squares technique can be used. The latter requires computing the partials of the loss function to be taken with respect to all of the unknowns. All partial derivatives for each segment and each unknown are provided below:
\begin{equation*}
\begin{aligned}
    \frac{\partial \p{s}{{\cal L}}_i}{\partial \p{s}{\B{\xi}_i}} &= \p{s}{\left(\ddot{\B{h}} - \ddot{\Omega}_1\B{h}_0 - \ddot{\Omega}_2\B{h}_f - \ddot{\Omega}_3\bdot{h}_0 - \ddot{\Omega}_4\bdot{h}_f \right)\T} \\ 
    \frac{\partial \p{1}{{\cal L}}_i}{\partial r_{1_i}} &= \p{1}{\ddot{\Omega}}_2 \\
    \frac{\partial \p{1}{{\cal L}}_i}{\partial v_{1_i}} &= \p{1}{\ddot{\Omega}}_4 \\  \frac{\partial \p{2}{{\cal L}}_i}{\partial r_{1_i}} &= \p{2}{\ddot{\Omega}}_1, \qquad  \frac{\partial \p{2}{{\cal L}}_i}{\partial v_{1_i}} = \p{2}{\ddot{\Omega}}_3 \\ 
    \frac{\partial \p{2}{{\cal L}}_i}{\partial r_{2_i}} &= \p{2}{\ddot{\Omega}}_2, \qquad \frac{\partial \p{2}{{\cal L}}_i}{\partial v_{2_i}} = \p{2}{\ddot{\Omega}}_4 \\
    \frac{\partial \p{3}{{\cal L}}_i}{\partial r_{2_i}} &= \p{3}{\ddot{\Omega}}_1 \\ \frac{\partial \p{3}{{\cal L}}_i}{\partial v_{2_i}} &= \p{3}{\ddot{\Omega}}_3
\end{aligned}
\end{equation*}
For the costate portion, if $i = j$
\begin{equation*}
   \frac{\partial {\cal L}_i}{\partial \B{\xi_{\lambda_i}}} = \beta(t) \left[ \left(\displaystyle\sum_{j=1}^3 \lambda^2_{v_j}\right)^{-1/2} - \lambda^2_{v_i} \, \left(\displaystyle\sum_{j=1}^3 \lambda^2_{v_j}\right)^{-3/2}\right] \B{h}_{\lambda}\T 
\end{equation*}
if $i \neq j$
\begin{equation*}
    \frac{\partial {\cal L}_i}{\partial \B{\xi_{\lambda_j}}} = \beta(t) \left[ - \lambda_{v_i} \, \lambda_{v_j} \, \left(\displaystyle\sum_{j=1}^3 \lambda^2_{v_j}\right)^{-3/2}\right] \B{h}_{\lambda}\T 
\end{equation*}
In addition to the loss functions for the problem dynamics given by Eq. (\ref{eq:v_dot}), a loss function associated with the transversality conditions for the Hamiltonian needs to be defined,
\begin{equation}
    \mathcal{L}_H = \alpha T_\text{max} + \sum_{i=1}^{3}\lambda_{v_i} a_{g_i} - \beta(t_f) \left(\sum_{i=1}^{3}\lambda^2_{v_i}\right)^{\frac{1}{2}}
\end{equation}
It follows that the only none-zero partial of $\mathcal{L}_H$ is with respect to $\B{\xi}_\lambda$, that is defined by
\begin{equation}
   \frac{\partial \mathcal{L}_H}{\partial \B{\xi_{\lambda_i}}} = \left[a_{g_i} - \beta(t_f) \, \lambda_{v_i} \, \left(\sum_{j=1}^3 \lambda^2_{v_j}\right)^{-1/2}\right] \B{h}_{\lambda}\T
\end{equation}
Next, by discretizing the domain over $N$ points, these partials become a vector or matrix where the second dimension is the number of unknowns. All partials can be combined into one augmented matrix and one augment vector for the loss function such that according to,
\begin{equation}\label{eq:aug_J}
    \mathbb{J} = \begin{bmatrix} \p{1}{J}_{\B{\xi}} & \p{1}{J}_{r_1, v_1}  & \B{0}_{(3N\times 3m)} & \B{0}_{(3N\times 6)} &  \B{0}_{(3N\times 3m)} &  \p{1}{J}_{\B{\xi}_\lambda} \\ \B{0}_{(3N\times 3m} & \p{2}{J}_{r_1, v_1} &   \p{2}{J}_{\B{\xi}} & \p{2}{J}_{r_2, v_2}  & \B{0}_{(3N\times 3m)} &  \p{2}{J}_{\B{\xi}_\lambda} \\ \B{0}_{(3N\times 3m)} & \B{0}_{(3N\times 6)} &  \B{0}_{(3N\times 3m)} & \p{3}{J}_{r_2, v_2} &  \p{3}{J}_{\B{\xi}}  & \p{3}{J}_{\B{\xi}_\lambda} \\ \B{0}_{(1\times 3m)} & \B{0}_{(1\times 6)} & \B{0}_{(1\times 3m)} & \B{0}_{(1\times 6)} & \B{0}_{(1\times 3m)} & J_H \end{bmatrix}_{(\{9N+1\}\times \{9m + 18\})}
\end{equation}
with the augmented vector of the loss functions and unknown vector defined as,
\begin{equation*}
\begin{aligned}
    \mathbb{L} &= \begin{Bmatrix} \p{1}{{\cal L}}_1\T & \p{1}{{\cal L}}_2\T & \p{1}{{\cal L}}_3\T & \p{2}{{\cal L}}_1\T & \p{2}{{\cal L}}_2\T & \p{2}{{\cal L}}_3\T & \p{3}{{\cal L}}_1\T & \p{3}{{\cal L}}_2\T & \p{3}{{\cal L}}_3\T & \mathcal{L}_H\end{Bmatrix}_{(\{9N+1\}\times 1)}\T \\ \Xi &= \begin{Bmatrix} \p{1}{\B{\xi}}_1\T & \p{1}{\B{\xi}}_2\T & \p{1}{\B{\xi}}_3\T & \B{r}_1\T & \B{v}_1\T & \p{2}{\B{\xi}}_1\T & \p{2}{\B{\xi}}_2\T & \p{2}{\B{\xi}}_3\T & \B{r}_2\T & \B{v}_2\T & \p{3}{\B{\xi}}_1\T & \p{3}{\B{\xi}}_2\T & \p{3}{\B{\xi}}_3\T & \B{\xi}_{\lambda_1}\T & \B{\xi}_{\lambda_2}\T & \B{\xi}_{\lambda_3}\T \end{Bmatrix}_{(9m + 18)}\T
\end{aligned}
\end{equation*}
The terms of Eq. (\ref{eq:aug_J}) are defined by the following equations:
\begin{equation*}
    \p{s}{J}_{\B{\xi}} = \begin{bmatrix} \frac{\partial \p{s}{{\cal L}}_1}{\partial \p{s}{\B{\xi}_1}} & \B{0} & \B{0} \\ \B{0} & \frac{\partial \p{s}{{\cal L}}_2}{\partial \p{s}{\B{\xi}_2}} & \B{0} \\ \B{0} & \B{0} & \frac{\partial \p{s}{{\cal L}}_3}{\partial \p{s}{\B{\xi}_3}} \end{bmatrix}_{(3N\times 3m)},
    \quad  \p{s}{J}_{\B{\xi}_\lambda} = \begin{bmatrix} J_{\B{\xi}_{\lambda_{11}}} & J_{\B{\xi}_{\lambda_{12}}} & J_{\B{\xi}_{\lambda_{13}}} \\ J_{\B{\xi}_{\lambda_{21}}} & J_{\B{\xi}_{\lambda_{22}}} & J_{\B{\xi}_{\lambda_{23}}} \\ J_{\B{\xi}_{\lambda_{31}}} & J_{\B{\xi}_{\lambda_{32}}} & J_{\B{\xi}_{\lambda_{33}}}  \end{bmatrix}_{(3N\times 6)}
\end{equation*}

\begin{equation*}
\begin{aligned}
    \p{1}{J}_{r_1,v_1} = \begin{bmatrix} \p{1}{ \ddot{\Omega}}_2 & \B{0} & \B{0} & \p{1}{ \ddot{\Omega}}_4 & \B{0} & \B{0} \\ \B{0} & \p{1}{ \ddot{\Omega}}_2 & \B{0} & \B{0} & \p{1}{ \ddot{\Omega}}_4 & \B{0} \\ \B{0} & \B{0} & \p{1}{ \ddot{\Omega}}_2 & \B{0} & \B{0} & \p{1}{ \ddot{\Omega}}_4 \end{bmatrix}_{(3N\times 6)} \\ \p{2}{J}_{r_1,v_1} = \begin{bmatrix} \p{1}{ \ddot{\Omega}}_1 & \B{0} & \B{0} & \p{1}{ \ddot{\Omega}}_3 & \B{0} & \B{0} \\ \B{0} & \p{1}{ \ddot{\Omega}}_1 & \B{0} & \B{0} & \p{1}{ \ddot{\Omega}}_3 & \B{0} \\ \B{0} & \B{0} & \p{1}{ \ddot{\Omega}}_1 & \B{0} & \B{0} & \p{1}{ \ddot{\Omega}}_3 \end{bmatrix}_{(3N\times 6)} \\
    \p{2}{J}_{r_2,v_2} = \begin{bmatrix} \p{2}{ \ddot{\Omega}}_2 & \B{0} & \B{0} & \p{2}{ \ddot{\Omega}}_4 & \B{0} & \B{0} \\ \B{0} & \p{2}{ \ddot{\Omega}}_2 & \B{0} & \B{0} & \p{2}{ \ddot{\Omega}}_4 & \B{0} \\ \B{0} & \B{0} & \p{2}{ \ddot{\Omega}}_2 & \B{0} & \B{0} & \p{2}{ \ddot{\Omega}}_4 \end{bmatrix}_{(3N\times 6)} \\
    \p{3}{J}_{r_2,v_2} = \begin{bmatrix} \p{3}{ \ddot{\Omega}}_1 & \B{0} & \B{0} & \p{3}{ \ddot{\Omega}}_3 & \B{0} & \B{0} \\ \B{0} & \p{3}{ \ddot{\Omega}}_1 & \B{0} & \B{0} & \p{3}{ \ddot{\Omega}}_3 & \B{0} \\ \B{0} & \B{0} & \p{3}{ \ddot{\Omega}}_1 & \B{0} & \B{0} & \p{3}{ \ddot{\Omega}}_3 \end{bmatrix}_{(3N\times 6)}
\end{aligned}
\end{equation*}
\begin{equation*}
    J_H = \begin{bmatrix} \dfrac{\partial \mathcal{L}_H}{\partial \B{\xi_{\lambda_1}}}, & \quad \dfrac{\partial \mathcal{L}_H}{\partial \B{\xi_{\lambda_2}}}, & \quad \dfrac{\partial \mathcal{L}_H}{\partial \B{\xi_{\lambda_3}}}\end{bmatrix}_{(1\times 6)}
\end{equation*}

Finally, using Eq (\ref{eq:aug_J}) along with the augment loss functions and unknown vector, an iterative least-square approach can be used to update the unknown parameters according to,
\begin{equation*}
    \Xi_{k+1}= \Xi_k - \dd \Xi_k \quad \text{where} \quad \dd \Xi_k  = \text{qr}(\mathbb{J}_k,\mathbb{L}_k)
\end{equation*}
Importantly, an initial estimate of the parameters is needed in order to initialize the iterative least-squares process. Since the problem is a boundary-value problem, a first guess for $\p{s}{\B{\xi}}_i, \B{r}_1, \B{r}_2, \B{v}_1$, and $\B{v}_2$ can be determined by simply connecting the initial and final position with a straight line and using this trajectory for a least squares fitting of the \ces\ describing the $\p{s}{r}_i$ terms. Next, since the $\lambda_{v_i}$ is related to the thrust direction it can be assumed similar to Eq. (51) of Ref. \cite{upg} such that,
\begin{equation*}
    \B{\lambda}_{v_0} = \frac{\B{v}_0}{||\B{v}_0||} 
\end{equation*}
However, the initialization of $\B{\lambda}_r = \B{0}$ is prohibitive in the sense of the TFC method because this involves setting $\B{\xi}_\lambda$ coefficients to zeros. Therefore in this paper the coefficients are initialized using the following prescription:
\begin{equation*}
    \B{\lambda}_{v_f} = -\frac{\B{r}_0}{||\B{r}_0||}
\end{equation*}

\section{Summary of Algorithm}
\begin{figure}[H]
    \centering\includegraphics[width=\linewidth]{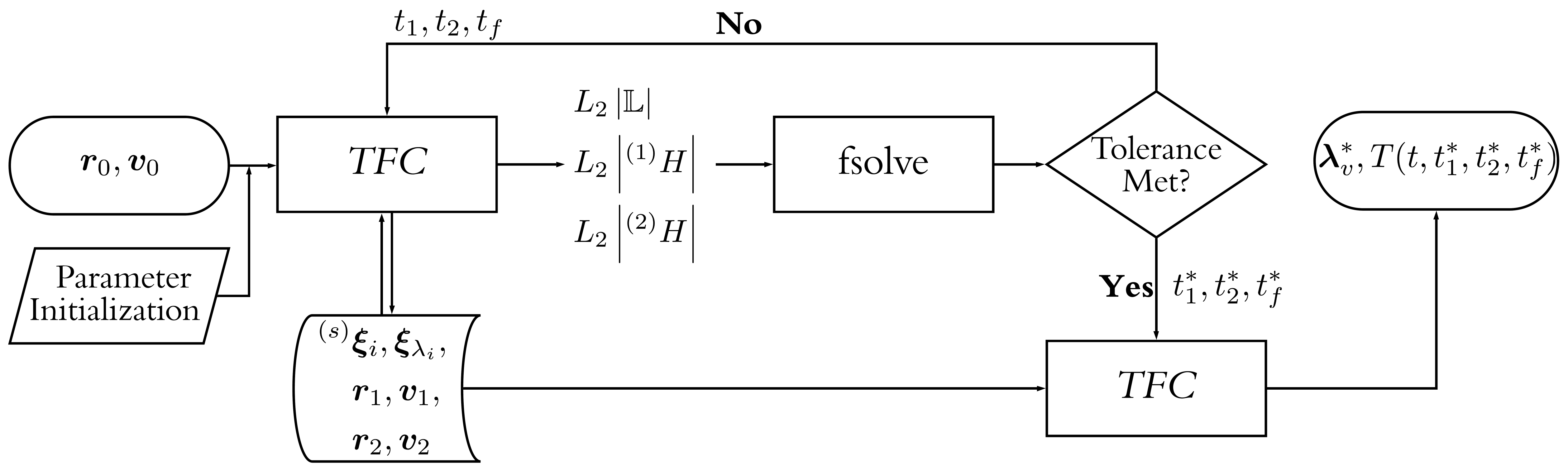}
    \caption{Summary of the full algorithm used with the TFC approach.}
    \label{fig:algorithm}
\end{figure}

Overall, the TFC method was used as a inner-loop function to minimize the residuals of the first-order necessary conditions subject to a prescribed thrust profile $T(t;t_1,t_2)$, i.e. the switching times  $t_1 $ and $t_2$ are assumed to be known by the TFC-based inner-loop routine. Consequently, an outer-loop routine has been developed to optimize the switching times given the $L_2$-norms of the residual of the first-order conditions, and the Hamiltonian over the first two segments (here, the MATLAB \cite{MATLAB} \verb"fsolve" was used) . A flow chart of the relevant inputs and outputs is provided in Fig. \ref{fig:algorithm}. 

\section{Results}
The proposed method was validated using two specific tests cases based on selected initial conditions defining a powered descent guidance scenario for landing on Mars. Within such framework, in Section \ref{sec:test1} the algorithm is tested on initial conditions where the optimal trajectory is characterized by a min-max thrust profile. Furthermore, in Section \ref{sec:test2} the case where the optimal thrust profile is max-min-max is studied. In both cases our results were compared with the GPOPS-II solutions for the boundary conditions. The algorithm was fully implemented in MATLAB R$2019a$ and all test cases were solved on a MacBook Pro (2016) macOS Version 10.15, with a 3.3 GHz Dual-Core Intel\textsuperscript{\textregistered} Core\texttrademark \, i7 and with 16 GB of RAM.

\subsection{Constant Test Parameters} 
We consider the trajectory optimization problem for a spacecraft performing powered descent for pinpoint landing on Mars. The gravitational field is assumed constant as generally, the powered descent starts below $1.5$ km. For the numerical test,  the lander parameters have been assumed to be similar to the one presented in Ref. \cite{acikmese2007convex} and reported in Table \ref{tab:test_parameters}.

\begin{table}[h]
\centering
\begin{tabular}{|c||c|}
\hline
Variable & Value \\
\hline\hline
$\B{a}_g$ [m/s$^2$] & $\begin{Bmatrix} 0, & 0, & -3.7114\end{Bmatrix}\T$ \\\hline
$I_{sp}$ [s] & $225$ \\\hline
$g_0$ [m/s$^2$] & $9.807$\\\hline 
$\overline{T}$ [N] & $3,100$\\\hline 
$N_T$ [-] & $6$\\\hline 
$\phi_T$ [deg] & $27$ \\\hline 
\hline
\end{tabular}
\caption{Constant parameters used in test cases.}
\label{tab:test_parameters}
\end{table}

Thrust magnitude bounds and the $\alpha$ parameter are defined as follows (see also Eqs. (\ref{eq:Tmin}-\ref{eq:alf})):
\begin{eqnarray}
    T_\text{min} &= 0.3 \overline{T}N_T\cos\phi_T &\approx 4,971.81 \: \: \: \text{[N]}\label{eq:Tmin}\\
    T_\text{max} &= 0.8 \overline{T}N_T\cos\phi_T &\approx 13,258.18 \: \text{[N]}\label{eq:Tmax}
\end{eqnarray}
where $ \overline{T}$ is the maximum thrust for a single engine, $N_T$ is the number of thrusters in the lander, and $\phi_T$ the cant angle of the thrusters with respect to the lander, and
\begin{equation}\label{eq:alf}
    \alpha = \frac{1}{I_{\rm sp} \, g_0 \, \cos\phi_T} \approx 5.0863 \cdot 10^{-4} \: \text{[s/m]}
\end{equation}
and where $I_{\rm sp}$ is the specific impulse of the engines and $g_0$ is Earth's gravitational constant.

\subsection{Test 1: Min-Max Trajectory}\label{sec:test1}
\begin{figure}[h]
    \centering\includegraphics[width=.70\linewidth]{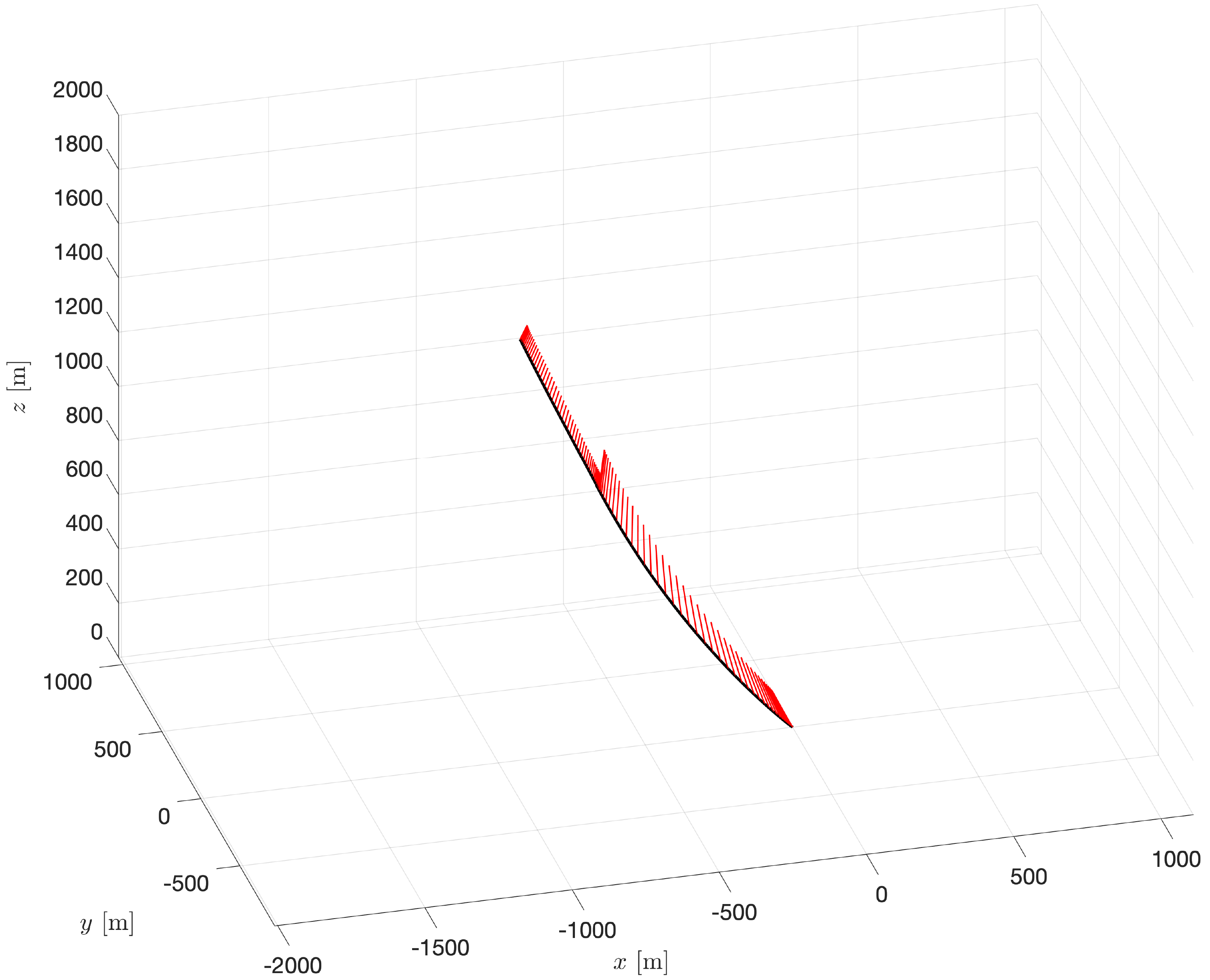}
    \caption{Landing trajectory for min-max thrust profile based on initial conditions, $\B{r}_0 = \protect\begin{Bmatrix} -900, & 10, & 1,500\protect\end{Bmatrix}\T $[m], $\B{v}_0 = \protect\begin{Bmatrix} 30, & -10, & -70\protect\end{Bmatrix}\T$ [m/s], $m_0 = 1,905$ [kg].}
    \label{fig:min_max_traj}
\end{figure}

For Test 1, initial conditions were selected such that the optimal thrust profile would be min-max (i.e. switch between minimum thrust to maximum thrust). Table \ref{tab:min_max} defines the boundary conditions for this test case and Fig. \ref{fig:min_max_traj} provides the converged trajectory using the TFC approach.
\begin{table}[h]
\centering
\begin{tabular}{|c||c|c|}
\hline
Variable & Initial & Final \\
\hline\hline
$\B{r}$ [m] & $\begin{Bmatrix} -900, & 10, & 1,500\end{Bmatrix}\T$ & $\begin{Bmatrix} 0, & 0, & 0\end{Bmatrix}\T$ \\\hline
$\B{v}$ [m/s] & $\begin{Bmatrix} 30, & -10, & -70\end{Bmatrix}\T$ & $\begin{Bmatrix} 0, & 0, & 0\end{Bmatrix}\T$ \\\hline
$m$ [kg] & $1,905$ & - \\\hline
\hline
\end{tabular}
\caption{Boundary conditions for min-max trajectory profile test case.}
\label{tab:min_max}
\end{table}

\begin{figure}[h]
    \centering\includegraphics[width=.9\linewidth]{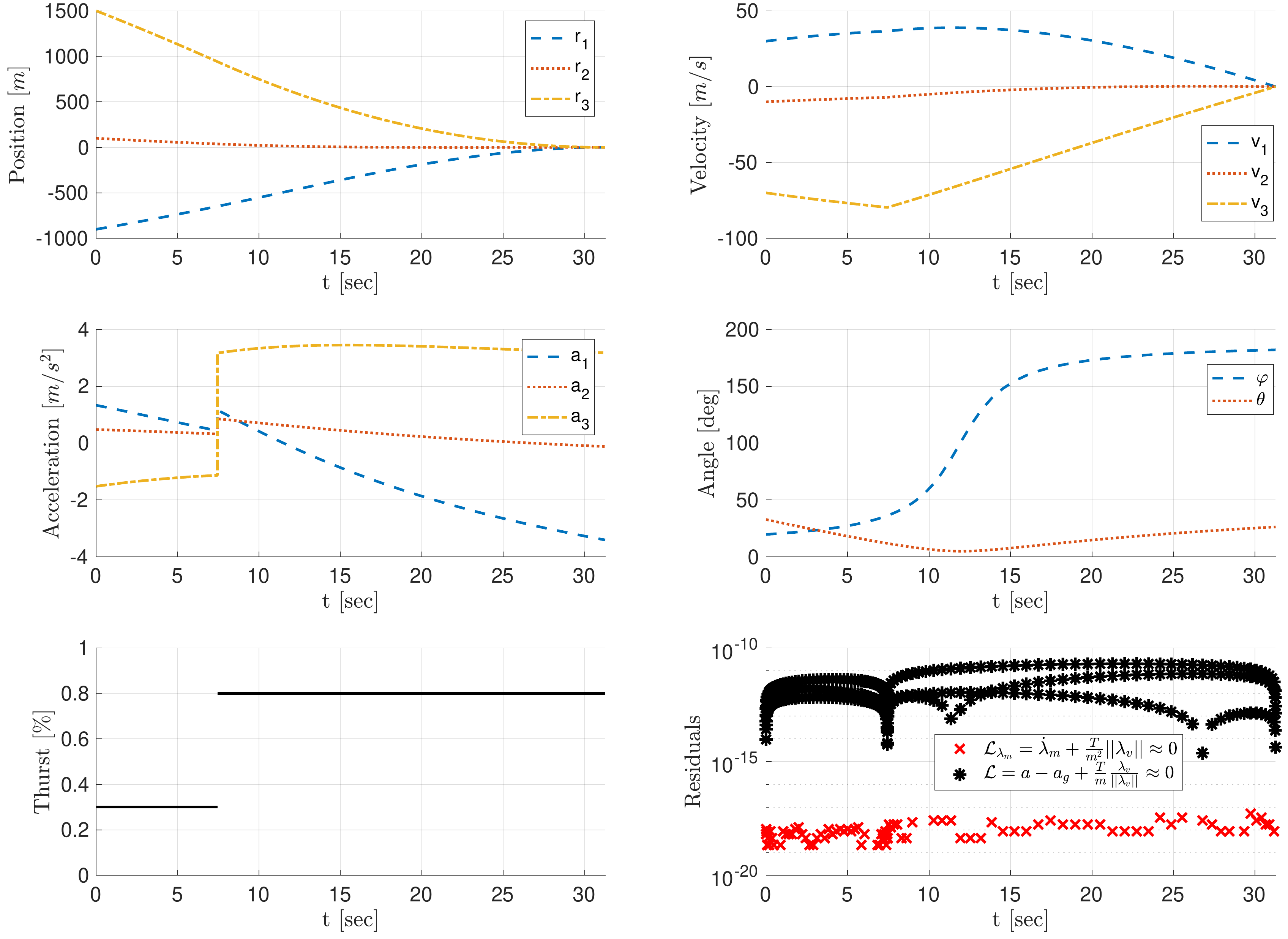}
    \caption{TFC solution to the min-max thrust profile case. The solution is presented in terms of the position, velocity, acceleration, and residuals of the differential equations.}
    \label{fig:min_max_tfc}
\end{figure}

\begin{figure}[h]
    \centering\includegraphics[width=.75\linewidth]{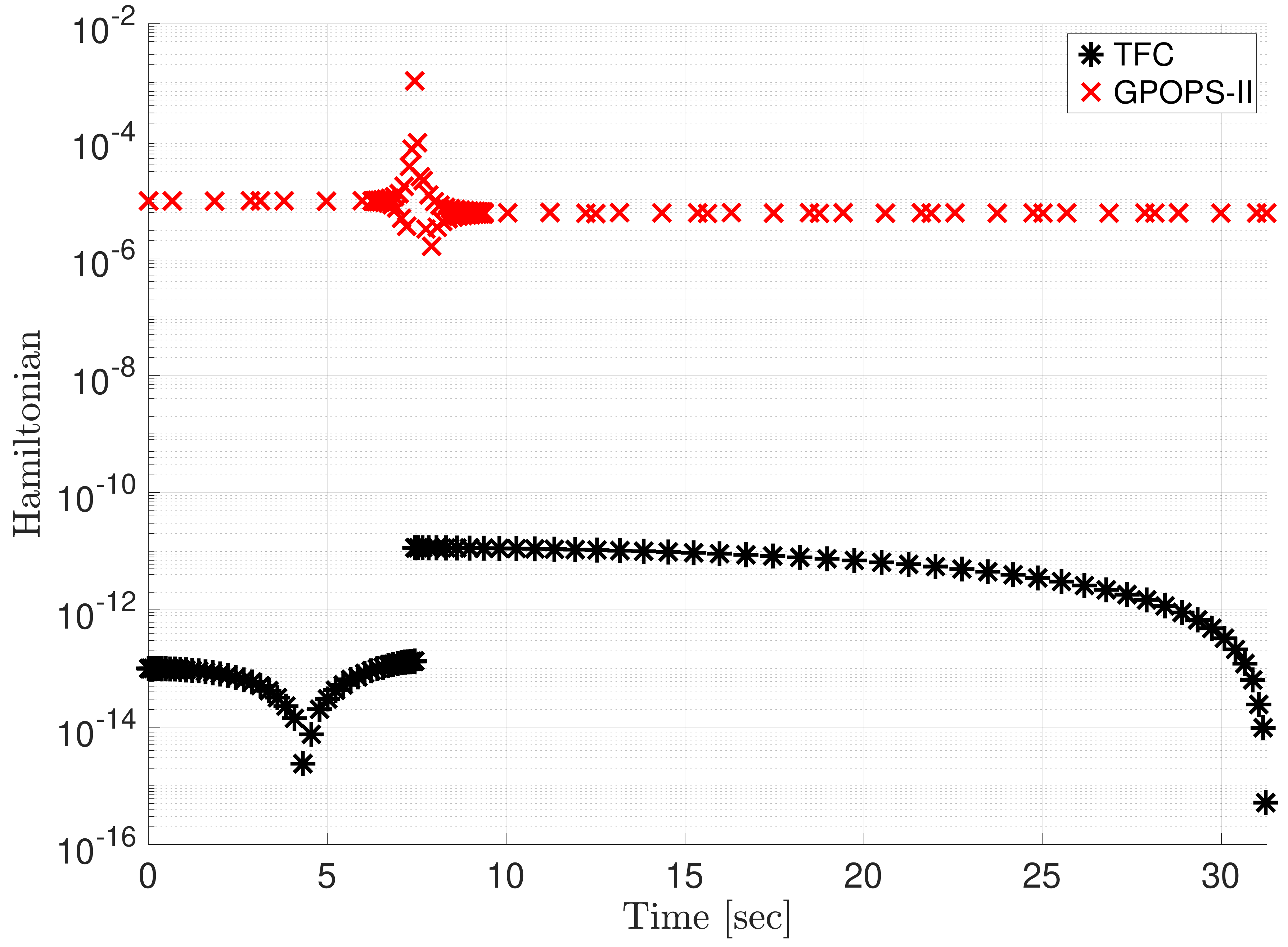}
    \caption{Comparison of Hamiltonian for TFC and GPOPS-II \cite{GPOPS} converged solutions for the min-max trajectory.}
    \label{fig:min_max_H}
\end{figure}

In addition to the trajectory, component plots of the position, velocity, and acceleration are provided in Fig. \ref{fig:min_max_tfc}. Furthermore, this figure also plots the residual of the governing differential equations for the mass and acceleration equation in order to quantify accuracy of the method. It can be seen that the TFC residual is about $\mathcal{O}(10^{-11})$ or less for the whole solution domain.
The accuracy of this approach was also compared to results obtained using GPOPS-II \cite{GPOPS} and is quantified in terms of the converged parameters, the $L_2$-norms of the Hamiltonian, and propellant mass used. Furthermore, to further justify accuracy of the solution, the converged parameters of initial costate values and switching times for each method were propagated using MATLAB \verb"ode45" a maximum tolerance in order to check the final position and velocity error and also the final error of the $\lambda_m$ term. The tabulated values of this test are provided in Table \ref{tab:min_max_results}. For this solution \verb"fsolve" iterated 27 times and over these iterations the TFC method converged in about 6 iterations. On average, the TFC innerloop took $76$ ms per call within the MATLAB implementation. Additionally, as a last point of comparison, the time history of the Hamiltonian for both methods is plotted versus time of the converged solution in Fig. \ref{fig:min_max_H}. 
\begin{table}[h]
\centering
\begin{tabular}{|c||c|c|c|}
\hline
Variable & TFC & GPOPS-II \cite{GPOPS} \\
\hline\hline
$L_2[\mathbb{L}]$ & $1.036 \cdot 10^{-10}$ & $-$  \\\hline
$L_2[H]$ & $5.488 \cdot 10^{-11}$ & $1.064 \cdot 10^{-3}$  \\\hline
$m_{\text{used}}$ [kg] & $179.447$  & $179.447$ \\\hline
$t_1$ [s] & $7.4430$ & $7.4430$ \\\hline 
$t_f$ [s] & $31.2623$ & $31.2623$  \\\hline 
$||\B{r}(t_f)||$ [m] & $2.886 \cdot 10^{-9}$ & $1.535\cdot 10^{-2}$  \\\hline 
$||\B{v}(t_f)||$ [m] & $3.166 \cdot 10^{-10}$ & $7.649\cdot 10^{-4}$  \\\hline 
$\lambda_m(t_f)$ [s] & $4.496 \cdot 10^{-14}$ & $-4.193\cdot 10^{-7}$  \\
\hline
\hline
\end{tabular}
\caption{Converged parameters for the TFC and GPOPS-II \cite{GPOPS} solution for the min-max trajectory test case.}
\label{tab:min_max_results}
\end{table}

\clearpage
\subsection{Test 2: Max-Min-Max Trajectory}\label{sec:test2}
\begin{figure}[h]
    \centering\includegraphics[width=.70\linewidth]{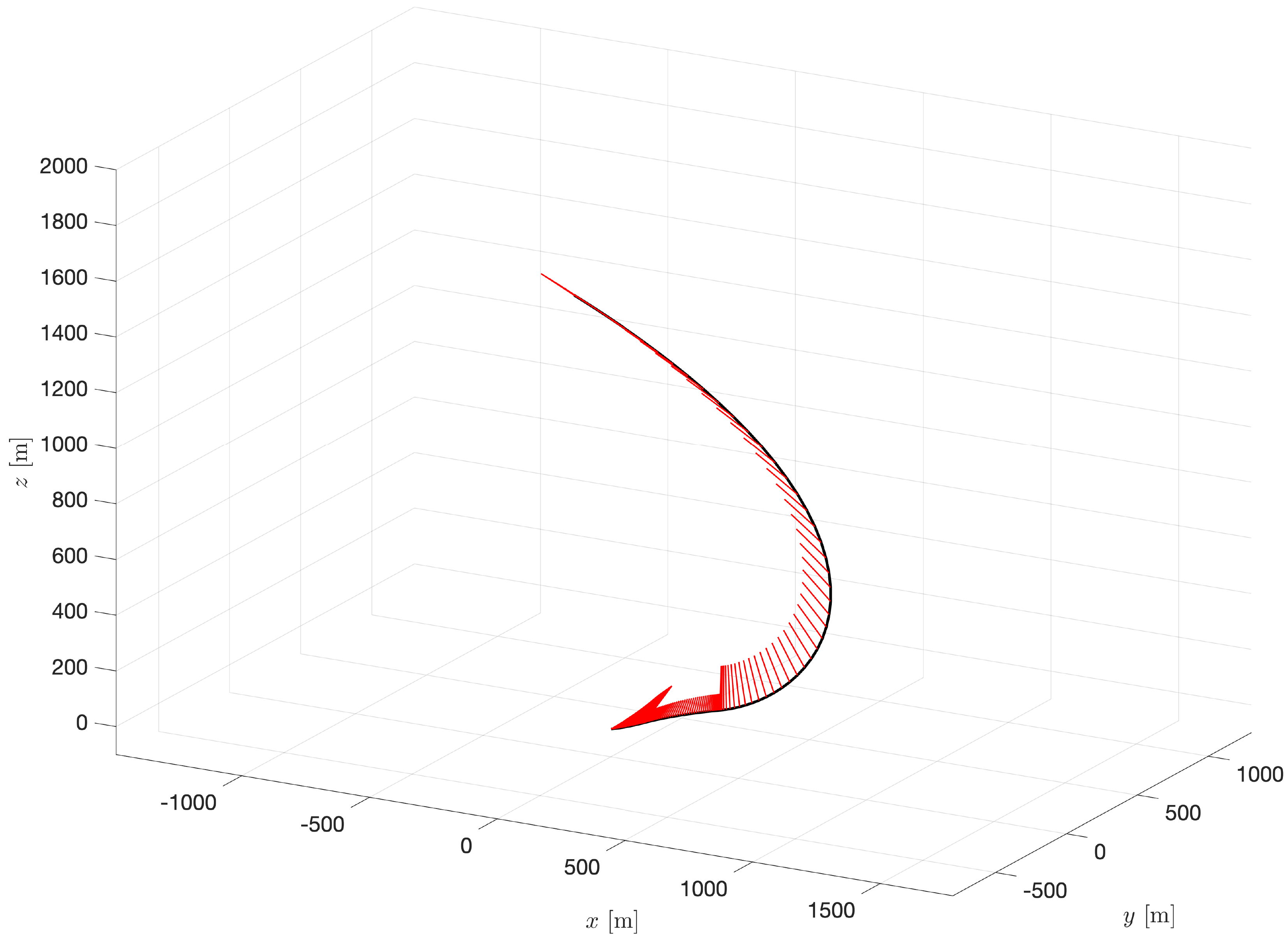}
    \caption{Landing trajectory for min-max thrust profile based on initial conditions, $\B{r}_0 = \protect\begin{Bmatrix} -200, & 100, & 1,500\protect\end{Bmatrix}\T $[m], $\B{v}_0 = \protect\begin{Bmatrix} 85, & -50, & -65\protect\end{Bmatrix}\T$ [m/s], $m_0 = 1,905$ [kg].}
    \label{fig:max_min_max_traj}
\end{figure}

\begin{figure}[h]
    \centering\includegraphics[width=.9\linewidth]{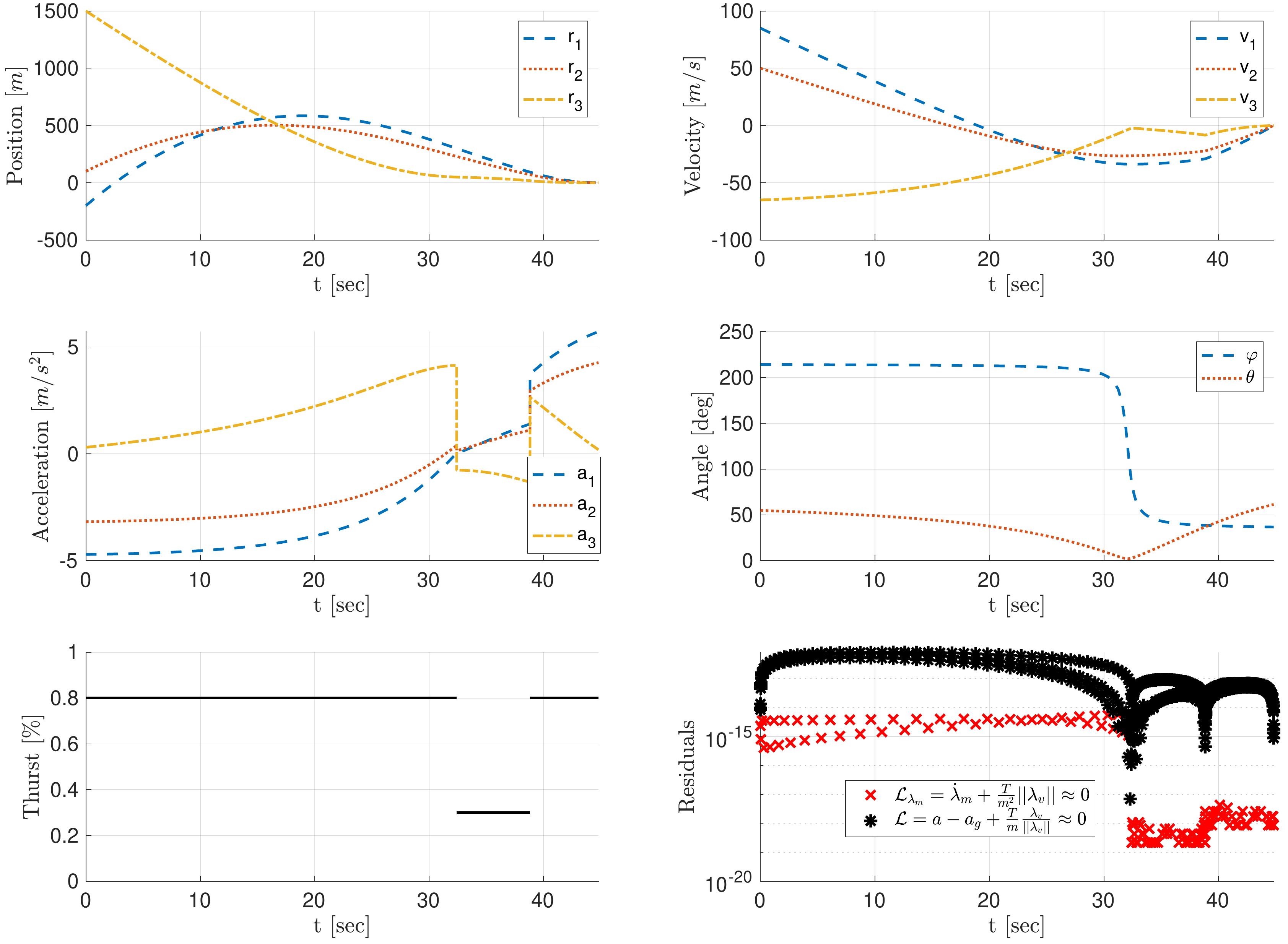}
    \caption{TFC solution to the min-max thrust profile case. The solution is presented in terms of the position, velocity, acceleration, and residuals of the differential equations.}
    \label{fig:max_min_max_tfc}
\end{figure}

\begin{figure}[h]
    \centering\includegraphics[width=.75\linewidth]{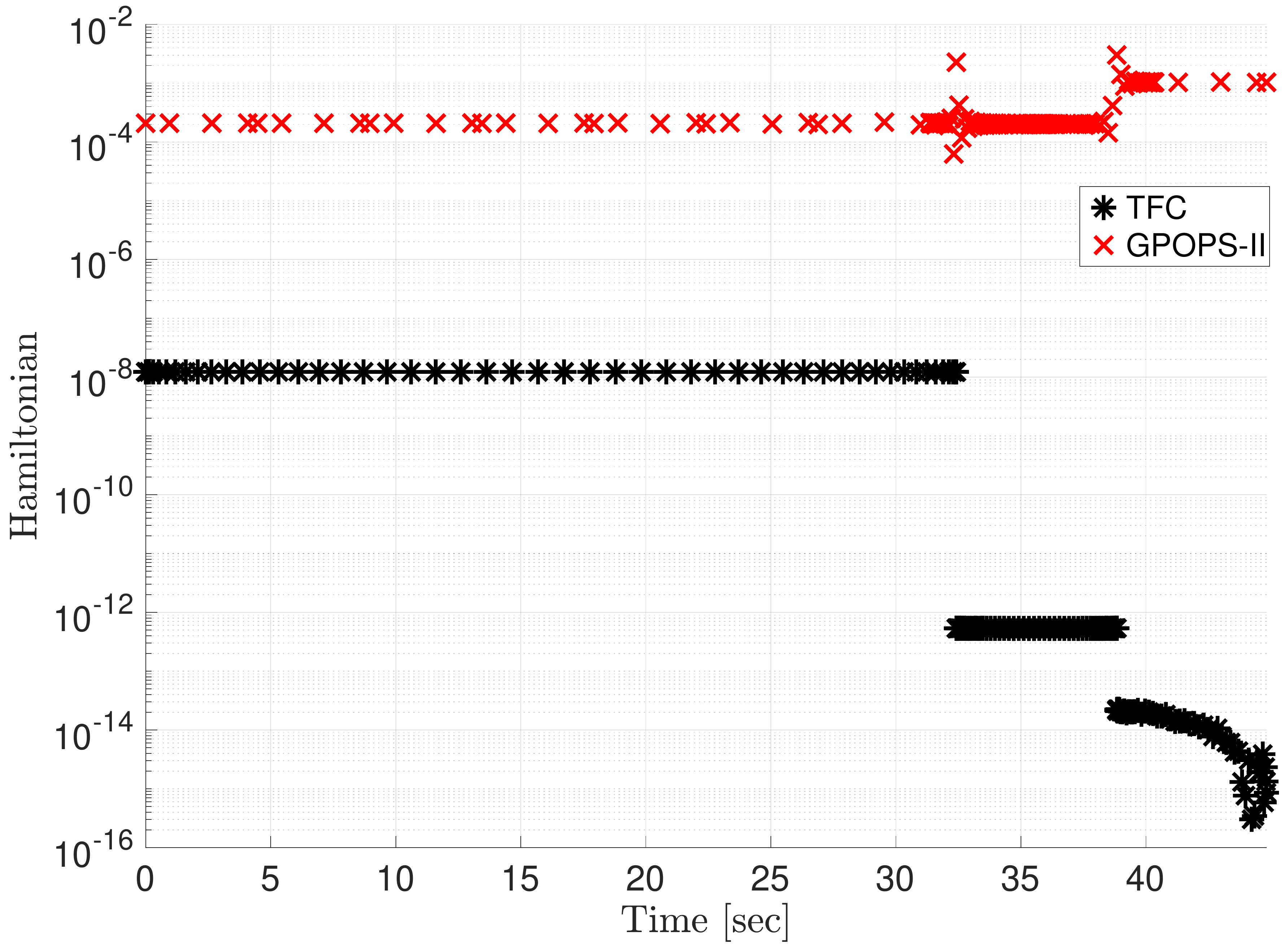}
    \caption{Comparison of Hamiltonian for TFC and GPOPS-II \cite{GPOPS} converged solutions for the max-min-max trajectory.}
    \label{fig:max_min_max_H}
\end{figure}

In test case 2, the initial conditions were specified such that the optimal solution exhibited a max-min-max profile (i.e. the thrust switches twice, max-min and min-max). The boundary conditions for this case are provided in Table \ref{tab:max_min_max} whereas Fig. \ref{fig:max_min_max_traj} reports the shape of the trajectory computed using the TFC-based algorithm.
\begin{table}[h]
\centering
\begin{tabular}{|c||c|c|}
\hline
Variable & Initial & Final \\
\hline\hline
$\B{r}$ [m] & $\begin{Bmatrix} -200, & 100, & 1,500\end{Bmatrix}\T$ & $\begin{Bmatrix} 0, & 0, & 0\end{Bmatrix}\T$ \\ \hline
$\B{v}$ [m/s] & $\begin{Bmatrix}85, & 50, & -65\end{Bmatrix}\T$ & $\begin{Bmatrix} 0, & 0, & 0\end{Bmatrix}\T$ \\ \hline
$m$ [kg] & $1,905$ & - \\ \hline
\hline
\end{tabular}
\caption{Boundary conditions for max-min-max trajectory profile test case.}
\label{tab:max_min_max}
\end{table}
Again, the TFC solution history is reported for each component of position, velocity, and acceleration in Fig. \ref{fig:max_min_max_tfc}. Additionally, the error is quantified by the residual of the governing equation of motion and the mass costate equation. It can be seen that the TFC residual is $\mathcal{O}(10^{-12})$ or less for the whole solution domain.
Similarly to the Test 1 case, the solution is compared with the one obtained via GPOPS-II \cite{GPOPS} for all converged parameters which now includes another switching time $t_2$. It can be seen that the magnitude of associated errors are similar to those presented in Section \ref{sec:test1}. For this solution \verb"fsolve" iterated 32 times whereas within the routine, the TFC algorithm converged in about 3 iterations. On average, the TFC inner loop took $81$ ms per call within the MATLAB implementation. Lastly, the Hamiltonian of the two methods is plotted as a function of time in Fig. \ref{fig:max_min_max_H} to highlight the optimality of both solutions.
\begin{table}[h]
\centering
\begin{tabular}{|c||c|c|c|}
\hline
Variable & TFC & GPOPS-II \cite{GPOPS} \\
\hline\hline
$L_2[\mathbb{L}]$ & $5.654\cdot 10^{-12}$ & $-$  \\\hline
$L_2[H]$ & $8.686\cdot 10^{-8}$ & $6.418 \cdot 10^{-3}$  \\\hline
$m_{\text{used}}$ [kg]& $275.205$ & $275.206$ \\\hline
$t_1$ [s] &$32.418$& $32.417$ \\\hline 
$t_2$ [s]&$38.838$ & $38.833$  \\\hline 
$t_f$ [s]&$44.823$ & $44.823$  \\\hline
$||\B{r}(t_f)||$ [m] & $8.330\cdot 10^{-10}$ & $1.350\cdot 10^{-1}$  \\\hline 
$||\B{v}(t_f)||$ [m] & $2.812\cdot 10^{-11}$ & $2.077\cdot 10^{-2}$  \\\hline 
$\lambda_m(t_f)$ [s] & $-8.815\cdot 10^{-15}$ & $-7.354\cdot 10^{-6}$  \\
\hline
\hline
\end{tabular}
\caption{Converged parameters for the TFC and GPOPS-II \cite{GPOPS} solution for the max-min-max trajectory test case.}
\label{tab:max_min_max_results}
\end{table}

\clearpage
\section{Conclusions}
In this paper, a new approach to solve the fuel-efficient powered descent guidance on large planetary bodies, like the Moon or Mars, via Theory of Functional Connections is proposed. The approach is based in solving the two-point boundary problem associated with the first-order necessary conditions which arise from the application of the indirect method. Through the formulation based on Theory of Functional Connections, the candidate solution is cast into a subspace that always satisfies the problem's constraints. 

For all test problems, the accuracy of the computed solutions are of $\mathcal{O}(10^{-10})$ or less. Additionally, since this method is based on simple numerical techniques (i.e. nonlinear least-squares), the proposed approach has the potential for solving the TPBVP arising from fuel-efficient solution in real-time. Indeed, each individual TFC-based, inner-loop solution is reported to converge within 6 iterations, each iterations ranging between 76 and 81 milliseconds. Accordingly, for all analyzed cases, the outer-loop converges in less than 32 iterations. Considering the overhead associated with the current MATLAB implementation, it is clear that timing performance can be substantially improved by using alternative coding languages.   Whereas the formulation presented in this work relies on an outer-loop optimizer to solve for the optimal thrust switching times (e.g. in this work, the MATLAB's \verb"fsolve" was employed), future works will be devoted to incorporate these parameters into a single loop. This can be done by incorporating the mapping coefficient $c$ as another unknown parameter that must be solved. This method has already been validated for solving free final time optimal control problems and it is expected to further increase the accuracy and speed of convergence of the algorithm.

Moreover, to reduce the sensitivity to the initial guess and further speed up the algorithm convergence, future works will investigate the possibility to use supervised machine learning techniques to learn proper initial guesses for the co-states, switching times, and final time as a function of the initial states.    

\section{Appendix}
\numberwithin{equation}{section}

In this paper, the solution of the mass costate equation was also solved using the TFC approach. Consider solving the equation,
\begin{equation}\label{eq:mass_costate}
    \dot{\lambda}_m = -\dfrac{T(t;t_1,t_2,t_f)}{m} ||\B{\lambda}_v|| \qquad \text{subject to:} \qquad \lambda_m(t_f) = 0,
\end{equation}
such that the thrust $T(t_f;t_1,t_2)$ case a switch in the dynamics of the problem at $t_1$ and $t_2$. Since this is simply a back propagation problem, each segment's initial condition is the final value of the prior section. Therefore, the \ces\ are defined as,
\begin{eqnarray}
    \p{1}{\lambda_m} &= (\B{h} - \B{h}_f)\T \B{\xi}_{\lambda_{m_1}} + \p{2}{\lambda_m}(t_1)\label{eq:cem_1}\\
    \p{2}{\lambda_m} &= (\B{h} - \B{h}_f)\T \B{\xi}_{\lambda_{m_2}} + \p{3}{\lambda_m}(t_2)\label{eq:cem_2}\\
    \p{3}{\lambda_m} &= (\B{h} - \B{h}_f)\T \B{\xi}_{\lambda_{m_3}} + \p{3}{\lambda_m}(t_f)\label{eq:cem_3}
\end{eqnarray}
Since all other parameters in Eq. (\ref{eq:mass_costate}) are known, the solution of, $\B{\xi}_{\lambda_{m_1}}, \B{\xi}_{\lambda_{m_2}}, \B{\xi}_{\lambda_{m_3}}$ can be solved by taking a derivative of Eqs. (\ref{eq:cem_1}-\ref{eq:cem_3}) and plugging them into Eq. (\ref{eq:mass_costate}) which produces three linear equations which can be discretized and must be solved backward from segment 3,
\begin{eqnarray}
    \bdot{h}\T \B{\xi}_{\lambda_{m_s}} = -\dfrac{T(t;t_1,t_2,t_f)}{m} ||\B{\lambda}_v|| \qquad \text{where} \qquad s = 1, 2, 3.
\end{eqnarray}
Although only a major equation were highlighted for the solution of the mass costate equation, the reader is directed to Refs. \cite{LDE,U-TFC,NDE} for more details on the TFC approach to solving ordinary differential equations and the advantages.

\bibliographystyle{plain} 
\bibliography{mybib} 

\end{document}